\documentclass[12pt]{article}%
\usepackage{graphicx}
\usepackage[intlimits]{amsmath}
\usepackage{latexsym}
\usepackage{amsfonts}
\usepackage{amssymb}%
 \makeatletter
\newfont{\footsc}{cmcsc10 at 8truept}
\newfont{\footbf}{cmbx10 at 8truept}
\newfont{\footrm}{cmr10 at 10truept}
\newtheorem{theorem}{Theorem}

\newtheorem{claim}{Claim}

\newtheorem{conjecture}[theorem]{Conjecture}

\newtheorem{lemma}[theorem]{Lemma}
\newtheorem{fact}[theorem]{Fact}

\newenvironment{proof}[1][Proof]{\noindent{\textbf {#1}  }}  {\hfill$\Box$\bigskip}

\begin{document}

\title{The maximum spectral radius of $C_{4}$-free graphs of given order and size}
\author{Vladimir Nikiforov\\{\small Department of Mathematical Sciences, University of Memphis, Memphis TN
38152}\\{\small email: vnikifrv@memphis.edu}}
\maketitle

\begin{abstract}
Suppose that $G$ is a graph with $n$ vertices and $m$ edges, and let $\mu$ be
the spectral radius of its adjacency matrix.

Recently we showed that if $G$ has no $4$-cycle, then $\mu^{2}-\mu\leq n-1,$
with equality if and only if $G$ is the friendship graph.

Here we prove that if $m\geq9$ and $G$ has no $4$-cycle, then $\mu^{2}\leq m,$
with equality if $G$ is a star. For $4\leq m\leq8$ this assertion
fails.\medskip

\textbf{Keywords: }\textit{4-cycles; graph spectral radius; graphs with no
4-cycles; friendship graph. }

\textbf{AMS classification: }05C50, 05C35.

\end{abstract}

This note is part of an ongoing project aiming to build extremal graph theory
on spectral grounds, see, e.g., \cite{BoNi07} and \cite{Nik02,Nik07h}.

Suppose $G$ is a graph with $n$ vertices and $m$ edges and let $\mu\left(
G\right)  $ be the spectral radius of its adjacency matrix. How large can
$\mu\left(  G\right)  $ be if $G$ has no cycles of length $4?$ This question
was partially answered in \cite{Nik07d}, Theorem 3:\medskip

\emph{Let }$G$\emph{ be a graph of order }$n$\emph{ with }$\mu\left(
G\right)  =\mu$\emph{. If }$G$\emph{ has no }$4$\emph{-cycles}$,$\emph{ then }%
\begin{equation}
\mu^{2}-\mu\leq n-1. \label{in3}%
\end{equation}
\emph{Equality holds if and only if every two vertices of }$G$\emph{ have
exactly one common neighbor.}\medskip

The condition for equality in (\ref{in3}) is a popular topic: as shown in
\cite{ERS66} and \cite{Hun02}, the only graph satisfying this condition is the
\emph{friendship graph} - a set of $\left\lfloor n/2\right\rfloor $ triangles
sharing a single common vertex. Thus equality is possible only for $n$ odd,
and (\ref{in3}) may be improved for even $n.$

\begin{conjecture}
Let $G$ be a graph of even order $n$ with $\mu\left(  G\right)  =\mu$. If $G$
has no $4$-cycles, then\emph{ }%
\begin{equation}
\mu^{3}-\mu^{2}-\left(  n-1\right)  \mu+1\leq0. \label{in4}%
\end{equation}
Equality holds if and only if $G$ is a star of order $n$ with $n/2-1$ disjoint
additional edges.
\end{conjecture}

Note that the number of edges of $G$ is missing in (\ref{in3}) and
(\ref{in4}). In contrast, Nosal \cite{Nos70} showed that if $\mu\left(
G\right)  >\sqrt{m},$ then $G$ has triangles. Our main result here is a
similar assertion for $4$-cycles:

\begin{theorem}
\label{th1}Let $m\geq9$ and $G$ be a graph with $m$ edges. If $\mu\left(
G\right)  >\sqrt{m},$ then $G$ has a $4$-cycle.
\end{theorem}

Note that Theorem \ref{th1} is tight, for all stars are $C_{4}$-free graphs
with $\mu\left(  G\right)  =\sqrt{m}.$ Also, let $S_{n,1}$ be the star of
order $n$ with an edge within its independent set: $S_{n,1}$ is $C_{4}$-free
and has $n$ edges, but $\mu\left(  G\right)  >\sqrt{n}$ for $4\leq n\leq8,$ as
shown in Lemma \ref{le2} below.

Observe that the original result of Nosal was sharpened in \cite{Nik07f},
Theorem 2, (i):\medskip

\emph{If }$\mu\left(  G\right)  \geq\sqrt{m},$\emph{ then }$G$\emph{ has a
triangle, unless }$G$\emph{ is a complete bipartite graph with possibly some
isolated vertices.}\medskip

It turns out that Theorem \ref{th1} can be sharpened likewise, at the price of
a considerably longer proof, which we omit.

\begin{theorem}
\label{th2}Let $m\geq9$ and $G$ be a graph with $m$ edges. If $\mu\left(
G\right)  \geq\sqrt{m},$ then $G$ has a $4$-cycle unless $G$ is a star or
$S_{9,1}$ with possibly some isolated vertices.
\end{theorem}

\subsection*{Proofs}

Our notation follows \cite{Bol98}; thus, if $G$ is a graph $G,$ and $X$ and
$Y$ are disjoint sets of vertices of $G$, we write:

- $E\left(  G\right)  $ for the edge set of $G$ and $e\left(  G\right)  $ for
$\left\vert E\left(  G\right)  \right\vert ;$

- $G\left[  X\right]  $ for the graph induced by $X,$ $E\left(  X\right)  $
for $E\left(  G\left[  X\right]  \right)  ,$ and $e\left(  X\right)  $ for
$\left\vert E\left(  X\right)  \right\vert ;$

- $e\left(  X,Y\right)  $ for the number of edges joining vertices in $X$ to
vertices in $Y;$

- $G-uv$ for the graph obtained by removing the edge $uv\in E\left(  G\right)
;$

- $\Gamma_{G}\left(  u\right)  $ for the set of neighbors of a vertex $u$ and
$d_{G}\left(  u\right)  $ for $\left\vert \Gamma_{G}\left(  u\right)
\right\vert ;$

- $\Gamma_{X}\left(  u\right)  $ for $\Gamma_{G}\left(  u\right)  \cap X$ and
$d_{X}\left(  u\right)  $ for $\left\vert \Gamma_{X}\left(  u\right)
\right\vert .$\medskip

We drop the subscript in $\Gamma_{G}\left(  u\right)  $ and $d_{G}\left(
u\right)  $ when it is understood.\medskip

Define $S_{n,k}$ to be the star of order $n$ with $k$ disjoint edges within
its independent set.\medskip

Next we give some facts, needed in the proof of Theorem \ref{th1}.

First, a fact implied by Theorem 1 in \cite{PaRe00}:

\begin{fact}
\label{maxen}Let $\mathbf{x}$ be a unit eigenvector to the spectral radius of
a graph with some edges. Then the entries of $\mathbf{x}$ do not exceed
$2^{-1/2}.\hfill\square$
\end{fact}

Next, a known fact, proved here for completeness:

\begin{lemma}
\label{le1}Let $A$ and $A^{\prime}$ be the adjacency matrices of two graphs
$G$ and $G^{\prime}$ on the same vertex set. Suppose that $\Gamma_{G}\left(
u\right)  \subsetneqq\Gamma_{G^{\prime}}\left(  u\right)  $ for some vertex
$u.$ If some positive eigenvector $\mathbf{x}$ to $\mu\left(  G\right)  $
satisfies $\left\langle A^{\prime}\mathbf{x},\mathbf{x}\right\rangle
\geq\left\langle A\mathbf{x},\mathbf{x}\right\rangle ,$ then $\mu\left(
G^{\prime}\right)  >\mu\left(  G\right)  .$
\end{lemma}

\begin{proof}
Since $\left\langle A^{\prime}\mathbf{x},\mathbf{x}\right\rangle
\geq\left\langle A\mathbf{x},\mathbf{x}\right\rangle ,$ the Rayleigh principle
implies that $\mu\left(  G^{\prime}\right)  \geq\mu\left(  G\right)  .$ If
$\mu\left(  G^{\prime}\right)  =\mu\left(  G\right)  ,$ then $\left\langle
A^{\prime}\mathbf{x},\mathbf{x}\right\rangle =\left\langle A\mathbf{x}%
,\mathbf{x}\right\rangle ,$ and, again by the Rayleigh principle, $\mathbf{x}$
is an eigenvector to $\mu\left(  G^{\prime}\right)  .$ But this is impossible,
for
\[
\mu\left(  G^{\prime}\right)  x_{u}=\sum_{uv\in E\left(  G^{\prime}\right)
}x_{v}>\sum_{uv\in E\left(  G\right)  }x_{v}=\mu\left(  G\right)  x_{u}.
\]
We use above that $\Gamma_{G}\left(  u\right)  \subset\Gamma_{G^{\prime}%
}\left(  u\right)  ,$ but there is some $v\in\Gamma_{G^{\prime}}\left(
u\right)  $ such that $v\notin\Gamma_{G}\left(  u\right)  .$ This completes
the proof of Lemma \ref{le1}.
\end{proof}

Finally, some facts about $\mu\left(  S_{n,k}\right)  $:

\begin{lemma}
\label{le2}(a) $\mu\left(  S_{n,k}\right)  $ is the largest root of the
equation
\[
x^{3}-x^{2}-\left(  n-1\right)  x+n-1-2k=0;
\]

(b) $\mu\left(  S_{n,k}\right)  \leq\sqrt{n-1+k}$ for $n-1+k\geq9,$ and
$\mu\left(  S_{n,1}\right)  >\sqrt{n}$ for $4\leq n\leq8.$
\end{lemma}

\begin{proof}
Suppose that $1$ is the dominating vertex of $S_{n,k}$, and $\left\{
2,3\right\}  ,\ldots,\left\{  2k,2k+1\right\}  $ are its $k$ additional edges.
Set $\mu=\mu\left(  S_{n,k}\right)  $ and let $\left(  x_{1},\ldots
,x_{n}\right)  $ be an eigenvector to $\mu.$ By symmetry,%
\[
x_{2}=x_{3}=\cdots=x_{2k+1}\text{ \ \ and \ \ }x_{2k+2}=x_{2k+3}=\cdots
=x_{n}.
\]
Setting $x_{1}=x,$ $x_{2}=y,$ $x_{n}=z,$ we see that
\begin{align*}
\mu z  &  =x,\\
\mu y  &  =y+x,\\
\mu x  &  =2ky+\left(  n-2k-1\right)  z.
\end{align*}
Solving this system, we find that $\mu$ is a root of the equation
\[
x^{3}-x^{2}-\left(  n-1\right)  x+n-1-2k=0.
\]
If $\mu$ is not the largest root of this equation, then it has to be smaller
than
\[
x_{\min}=1/3+\sqrt{1/9+\left(  n-1\right)  /3},
\]
the point where the function
\[
f_{k}\left(  x\right)  =x^{3}-x^{2}-\left(  n-1\right)  x+n-1-2k
\]
has a local minimum. This, however, is not possible since
\[
\mu>\sqrt{n-1}>1/3+\sqrt{1/9+\left(  n-1\right)  /3}.
\]
This completes the proof of \emph{(a),}

To prove \emph{(b)} note that
\begin{align*}
f_{k}\left(  \sqrt{n-1+k}\right)   &  =\left(  \sqrt{n-1+k}\right)
^{3}-\left(  \sqrt{n-1+k}\right)  ^{2}-\left(  n-1\right)  \sqrt{n-1+k}+n-2k\\
&  =k\left(  \sqrt{n-1+k}-3\right)  ,
\end{align*}
implying the assertion since $\sqrt{n-1+k}>x_{\min}$ and $f_{k}\left(
x\right)  $ is increasing for $x>x_{\min}.$
\end{proof}

\subsubsection*{Proof of Theorem \textbf{\ref{th1}}}

Let $m\geq9,$ and assume for a contradiction that $G$ is a $C_{4}$-free graph
with $m$ edges, satisfying $\mu\left(  G\right)  >\sqrt{m}.$ Set $\mu
=\mu\left(  G\right)  ,$ and suppose that%
\begin{equation}
\mu=\max\left\{  \mu\left(  G\right)  :G\text{ is a }C_{4}\text{-free graph
with }e\left(  G\right)  =m\right\}  . \label{cond}%
\end{equation}
Also, for the purposes of the proof we may and shall suppose that $G$ has no
isolated vertices. This implies that $G$ is connected.

Indeed, let $G_{1}$ be a component of $G$ with $\mu\left(  G_{1}\right)
=\mu\left(  G\right)  ,$ and let $G_{2}$ be the nonempty union of the
remaining components of $G.$ Remove an edge from $G_{2},$ and add an edge
between $G_{1}$ and $G_{2}.$ The resulting graph is $C_{4}$-free with $m$
edges, but its spectral radius is larger than $\mu,$ contradicting
(\ref{cond}). Hence, $G$ is connected.

The essentially part of the proof is induction on $m,$ but it needs some
preparation. We first introduce some structure in $G$ and settle several cases
with direct arguments, in particular the case $m\leq13$. Then, having
restricted the structure of $G,$ we prove the induction step. Now the details.

Let $\left\{  1,\ldots,n\right\}  $ be the vertices of $G,$ and let
$\mathbf{x}=\left(  x_{1},\ldots,x_{n}\right)  $ be a positive unit
eigenvector to $\mu$, i.e.,
\[
\mu=2\sum_{ij\in E\left(  G\right)  }x_{i}x_{j}.
\]

By symmetry, suppose that $x_{1}\geq\cdots\geq x_{n}.$ We claim that all
vertices of degree $1$ are joined to vertex $1.$

Indeed, assume for a contradiction that there exists a vertex $u\neq1$ such
that $d\left(  u\right)  =1$ and $u$ is joined to $v\neq1.$ Remove the edge
$uv$ and join $u$ to vertex $1.$ The resulting graph $G^{\prime}$ is $C_{4}%
$-free and has $m$ edges. Also, we see that
\[
\sum_{ij\in E\left(  G^{\prime}\right)  }x_{i}x_{j}=\sum_{ij\in E\left(
G\right)  }x_{i}x_{j}+x_{u}\left(  x_{1}-x_{v}\right)  \geq\sum_{ij\in
E\left(  G\right)  }x_{i}x_{j}.
\]
Since $\Gamma_{G}\left(  1\right)  \subsetneqq\Gamma_{G^{\prime}}\left(
1\right)  ,$ Lemma \ref{le1} implies that $\mu\left(  G^{\prime}\right)
>\mu,$ contradicting (\ref{cond}). Hence, all vertices of degree $1$ are
joined to vertex $1.$

Let $A=\left(  a_{ij}\right)  $ be the adjacency matrix of $G$ and
$A^{2}=B=\left(  b_{ij}\right)  .$ Since $\mathbf{x}$ is an eigenvector of $B$
to $\mu^{2},$ we have
\begin{equation}
x_{1}\mu^{2}=\sum_{i=1}^{n}b_{1i}x_{i}\leq x_{1}\sum_{i=1}^{n}b_{1i}%
=\sum_{i=1}^{n}\sum_{j=1}^{n}a_{1j}a_{ji}=x_{1}\sum_{v\in\Gamma\left(
1\right)  }d\left(  v\right)  . \label{in2}%
\end{equation}

Set
\[
U=\Gamma\left(  1\right)  ,\text{ \ \ \ }W=\left\{  2,3,\ldots,n\right\}
\backslash\Gamma\left(  1\right)  ,
\]
and let $t=e\left(  U\right)  $ and $q=e\left(  W\right)  .$ We see that
\[
\sum_{v\in U}d\left(  v\right)  =d\left(  1\right)  +2e\left(  U\right)
+e\left(  U,W\right)  =e\left(  G\right)  -e\left(  W\right)  +e\left(
U\right)  =m-q+t.
\]
Thus (\ref{in2}) gives $\mu^{2}\leq m+t-q,$ and from $\mu^{2}>m,$ we get the
crucial inequality $t\geq q+1$.

Since all vertices of degree $1$ belong to $U,$ we have $d\left(  u\right)
\geq2$ for all $u\in W.$ Also, since $G$ is $C_{4}$-free, a vertex in $W$ can
be joined to at most one vertex in $U.$ Thus, for all $w\in W$ we have
$d_{W}\left(  w\right)  \geq d\left(  w\right)  -1\geq1,$ and consequently,
\[
2q=\sum_{w\in W}d_{W}\left(  w\right)  \geq\sum_{w\in W}1=\left\vert
W\right\vert .
\]

Suppose first that $q=0.$ Then $\left\vert W\right\vert =0,$ and so, $e\left(
U,W\right)  =0.$ Therefore, vertex $1$ is dominating and $G=S_{m+1-t,t}$. By
Lemma \ref{le2},%
\[
\mu=\mu\left(  S_{m+1-t,t}\right)  \leq\sqrt{m}%
\]
for $m\geq9,$ contradicting the hypothesis. Therefore, $q\geq1.$

The next claim gives a useful property of $G\left[  W\right]  ,$ and, in
particular, settles the case $q=1.$

\begin{claim}
\label{cl0}The graph $G\left[  W\right]  $ contains no isolated edges.
\end{claim}

\begin{proof}
Let $uv\in E\left(  W\right)  $ be an isolated edge. Since $d\left(  u\right)
\geq2$ and $d\left(  v\right)  \geq2,$ we see that $d_{U}\left(  u\right)
=d_{U}\left(  v\right)  =1.$ Let $\left\{  k\right\}  =\Gamma_{U}\left(
u\right)  $ and $\left\{  l\right\}  =\Gamma_{U}\left(  v\right)  .$ Remove
the edges $uk,vl,$ and join $u$ and $v$ to the vertex $1.$ The resulting graph
$G^{\prime}$ is $C_{4}$-free and has $m$ edges. Also, we see that
\[
\sum_{ij\in E\left(  G^{\prime}\right)  }x_{i}x_{j}=\sum_{ij\in E\left(
G\right)  }x_{i}x_{j}+x_{u}\left(  x_{1}-x_{k}\right)  +x_{v}\left(
x_{1}-x_{l}\right)  \geq\sum_{ij\in E\left(  G\right)  }x_{i}x_{j}.
\]
Since $\Gamma_{G}\left(  1\right)  \subsetneqq\Gamma_{G^{\prime}}\left(
1\right)  ,$ Lemma \ref{le1} implies that $\mu\left(  G^{\prime}\right)
>\mu,$ contradicting (\ref{cond}), and completing the proof of Claim \ref{cl0}.
\end{proof}

\bigskip

Claim \ref{cl0} implies that $q\geq2.$ Our next goal is to obtain a
contradiction for $m\leq13.$ Indeed, suppose that $m\leq13;$ then $q\geq2$
gives%
\[
13\geq m=3t+e\left(  U,W\right)  +q\geq4q+3+e\left(  U,W\right)
\geq11+e\left(  U,W\right)  ,
\]
which is possible only if $q=2,$ $e\left(  U,W\right)  \leq2,$ and $t=3.$

The graph $G\left[  W\right]  $ has $2$ non-isolated edges, and thus is a path
of order $3.$ Let $u,v,w$ be the vertices of this path and suppose that $uv\in
E\left(  W\right)  $ and $vw\in E\left(  W\right)  .$ Since $d\left(
u\right)  \geq2$ and $d\left(  w\right)  \geq2,$ we find that $d_{U}\left(
u\right)  =d_{U}\left(  w\right)  =1.$ This, in view of $e\left(  U,W\right)
\leq2,$ gives $e\left(  U,W\right)  =2,$ and so, $v$ has no neighbors in $U.$

Let $\left\{  k\right\}  =\Gamma_{U}\left(  u\right)  $ and $\left\{
l\right\}  =\Gamma_{U}\left(  w\right)  .$ Remove the edges $uk,wl,uv,$ and
join $u,v,w$ to the vertex $1.$ The resulting graph $G^{\prime}$ is $C_{4}%
$-free and has $m$ edges. Also, we see that
\[
\sum_{ij\in E\left(  G^{\prime}\right)  }x_{i}x_{j}=\sum_{ij\in E\left(
G\right)  }x_{i}x_{j}+x_{u}\left(  x_{1}-x_{k}\right)  +x_{w}\left(
x_{1}-x_{l}\right)  +x_{v}\left(  x_{1}-x_{u}\right)  \geq\sum_{ij\in E\left(
G\right)  }x_{i}x_{j}.
\]
Since $\Gamma_{G}\left(  1\right)  \nsubseteq\Gamma_{G^{\prime}}\left(
1\right)  ,$ Lemma \ref{le1} implies that $\mu\left(  G^{\prime}\right)
>\mu,$ contradicting (\ref{cond}).

At this point we have proved the theorem for $9\leq m\leq13.$ Assume now that
$m\geq14$ and that the theorem holds for $m-1;$ we shall prove it for $m$. The
induction step is based on three claims.

\begin{claim}
\label{cl1}If an edge $uv\in E\left(  G\right)  $ satisfies $d\left(
u\right)  =d\left(  v\right)  =2,$ then $x_{u}x_{v}<1/4\mu$.
\end{claim}

\begin{proof}
Let $\left\{  i,u\right\}  =\Gamma\left(  v\right)  $ and $\left\{
j,v\right\}  =\Gamma\left(  u\right)  .$ From
\[
\mu x_{u}=x_{i}+x_{v}\leq x_{1}+x_{v}\text{ \ \ and \ }\mu x_{v}\leq
x_{1}+x_{u}\leq x_{1}+x_{v}%
\]
we see that $x_{u}+x_{v}=2x_{1}/\left(  \mu-1\right)  .$ Hence, using the
AM-QM inequality and Fact \ref{maxen}, we obtain
\[
x_{u}x_{v}\leq\left(  \frac{x_{u}+x_{v}}{2}\right)  ^{2}=\frac{x_{1}^{2}%
}{\left(  \mu-1\right)  ^{2}}\leq\frac{1}{2\left(  \mu-1\right)  ^{2}}%
\leq\frac{1}{4\mu}%
\]
whenever $\mu^{2}\geq14.$ This completes the proof of Claim \ref{cl1}.
\end{proof}

\begin{claim}
\label{cl2}Let $m\geq20.$ Let the vertices $u,v,w$ satisfy $d\left(  u\right)
=d\left(  w\right)  =2$ and $d\left(  v\right)  =3,$ and let $v$ be joined to
$u$ and $w$. Then either $x_{u}x_{v}<1/4\mu$ or $x_{w}x_{v}<1/4\mu.$
\end{claim}

\begin{proof}
We first note that if $x\geq\sqrt{20},$ then
\begin{equation}
\frac{\left(  x^{2}-2\right)  ^{2}}{x\left(  x+1\right)  \left(  x+2\right)
}>\frac{x^{4}-4x^{2}}{x\left(  x+1\right)  \left(  x+2\right)  }%
=\frac{x\left(  x-2\right)  }{x+1}=\frac{x^{2}-4x-2}{x+1}+2>2. \label{in1}%
\end{equation}

Next, letting $\Gamma\left(  u\right)  =\left\{  i,v\right\}  $,
$\Gamma\left(  w\right)  =\left\{  j,v\right\}  ,$ and $\Gamma\left(
v\right)  =\left\{  k,u,w\right\}  ,$ we see that%
\begin{align*}
\mu x_{u}  &  =x_{i}+x_{v}\leq x_{1}+x_{v},\\
\mu x_{w}  &  =x_{j}+x_{v}\leq x_{1}+x_{v},\\
\mu x_{v}  &  =x_{k}+x_{u}+x_{w}\leq x_{1}+x_{u}+x_{w},
\end{align*}
and therefore,%
\begin{align*}
\mu\left(  x_{u}+x_{w}\right)   &  \leq x_{1}+2x_{v},\\
\mu x_{v}  &  \leq x_{1}+x_{u}+x_{w}.
\end{align*}
The solution of this system is%
\[
x_{u}+x_{w}\leq2\frac{\mu+1}{\mu^{2}-2}x_{1},\text{ \ \ }x_{v}\leq\frac{\mu
+2}{\mu^{2}-2}x_{1}.
\]
Now, assuming $x_{u}\geq x_{w},$ and using Fact \ref{maxen}, we obtain
\[
x_{u}x_{v}\leq\frac{\left(  \mu+1\right)  \left(  \mu+2\right)  }{\left(
\mu^{2}-2\right)  ^{2}}x_{1}^{2}\leq\frac{\left(  \mu+1\right)  \left(
\mu+2\right)  }{2\left(  \mu^{2}-2\right)  ^{2}}.
\]
Finally, inequality (\ref{in1}) implies that
\[
x_{u}x_{v}\leq\frac{\left(  \mu+1\right)  \left(  \mu+2\right)  }{2\left(
\mu^{2}-2\right)  ^{2}}\leq\frac{1}{4\mu}%
\]
whenever $\mu^{2}\geq20.$ This completes the proof of Claim \ref{cl2}.\medskip
\end{proof}

\begin{claim}
\label{cl3}If there exists $uv\in E\left(  G\right)  $ satisfying $x_{u}%
x_{v}\leq1/4\mu,$ then $\mu^{2}\left(  G-uv\right)  >\mu^{2}-1.$
\end{claim}

\begin{proof}
For every edge $uv\in E\left(  G\right)  ,$ by the Rayleigh principle, we have%
\[
\mu^{2}\left(  G-uv\right)  \geq\left(  2\sum_{ij\in E\left(  G-uv\right)
}x_{i}x_{j}\right)  ^{2}=\left(  \mu-2x_{u}x_{v}\right)  ^{2}>\mu^{2}-4\mu
x_{u}x_{v}\geq\mu^{2}-1,
\]
completing the proof of Claim \ref{cl3}.
\end{proof}

\bigskip

Having proved the claims, we proceed with the induction step. If there exists
$uv\in E\left(  U\right)  $ with $d\left(  u\right)  =d\left(  v\right)  =2,$
then by Claims \ref{cl1} and \ref{cl3} we obtain $\mu\left(  G-uv\right)
>\sqrt{m-1};$ by the induction hypothesis $G$ contains a $C_{4},.a$ contradiction.

Hereafter, we assume that $d\left(  u\right)  +d\left(  v\right)  \geq5$ for
all $uv\in E\left(  U\right)  .$ For every edge $uv\in E\left(  U\right)  ,$
let $W_{uv}=\Gamma_{W}\left(  u\right)  \cup\Gamma_{W}\left(  v\right)  .$
Since a vertex in $W$ can be joined to at most one vertex in $U,$ the sets
$W_{uv},$ $uv\in E\left(  U\right)  $ are disjoint. From
\[
2q=2e\left(  W\right)  =\sum_{w\in W}d_{W}\left(  w\right)  \geq\sum_{uv\in
E\left(  U\right)  }\sum_{w\in W_{uv}}d_{W}\left(  w\right)  \geq t\min_{uv\in
E\left(  U\right)  }\sum_{w\in W_{uv}}d_{W}\left(  w\right)
\]
we see that there is an edge $uv\in E\left(  U\right)  $ such that $\sum_{w\in
W_{uv}}d_{W}\left(  w\right)  \leq1.$ Then from
\[
\left\vert W_{uv}\right\vert =d\left(  u\right)  +d\left(  v\right)  -4\geq1
\]
we conclude that $W_{uv}$ contains a single vertex $w,$ and that $d_{W}\left(
w\right)  =1.$

Assume, by symmetry, that $w$ is joined to $v.$ Then, $d\left(  u\right)  =2,$
$d\left(  w\right)  =2,$ and $d\left(  v\right)  =3.$ Now, if $m\geq20,$
Claims \ref{cl2} and \ref{cl3} imply either $\mu\left(  G-vw\right)
>\sqrt{m-1}$ or $\mu\left(  G-uv\right)  >\sqrt{m-1};$ by the induction
hypothesis $G$ contains a $C_{4},$ contradiction.

To complete the proof we have to settle the case when $15\leq m\leq19$ and
$d\left(  u\right)  +d\left(  v\right)  \geq5$ holds for all $uv\in E\left(
U\right)  .$ We shall show that these conditions also lead to a contradiction.

From
\[
e\left(  U,W\right)  =\sum_{uv\in E\left(  U\right)  }d_{W}\left(  u\right)
+d_{W}\left(  v\right)  \geq\sum_{uv\in E\left(  U\right)  }\left(
5-4\right)  =t
\]
and
\begin{equation}
19\geq m=3t+e\left(  U,W\right)  +q\geq3t+t+q\geq5q+4 \label{in5}%
\end{equation}
we see that $q\leq3$ and $t\leq4.$

Consider first the case $q=3.$ From (\ref{in5}) we find that this is possible
only if $m=19,$ $t=4,$ $e\left(  U,W\right)  =4.$ This implies also that
$\left\vert W\right\vert \geq e\left(  U,W\right)  \geq4.$

$G\left[  W\right]  $ has no isolated vertices and, by Claim \ref{cl0}, it has
no isolated edges either. Thus, from $e\left(  W\right)  =3$ we see that
$G\left[  W\right]  $ is a tree of order $4.$ Now the structure of $G$ is
determined: $G$ consists of $4$ triangles sharing vertex $1$, a tree $T$ of
order $4,$ and a $4$-matching joining every vertex of $T$ to a separate triangle.

Select $u\in W$ to be with $d_{W}\left(  u\right)  =1$ and let $\left\{
v\right\}  =\Gamma_{W}\left(  u\right)  ,$ $\left\{  k\right\}  =\Gamma
_{U}\left(  u\right)  ,$ $\left\{  l\right\}  =\Gamma_{U}\left(  v\right)  .$
Suppose that $x_{k}\geq x_{l},$ remove the edge $vl,$ and add the edge $vk.$
The resulting graph $G^{\prime}$ is $C_{4}$-free and has $m$ edges. Also, we
see that
\[
\sum_{ij\in E\left(  G^{\prime}\right)  }x_{i}x_{j}=\sum_{ij\in E\left(
G\right)  }x_{i}x_{j}+x_{v}\left(  x_{k}-x_{l}\right)  \geq\sum_{ij\in
E\left(  G\right)  }x_{i}x_{j}.
\]
Since $\Gamma_{G}\left(  k\right)  \subsetneqq\Gamma_{G^{\prime}}\left(
k\right)  ,$ Lemma \ref{le1} implies that $\mu\left(  G^{\prime}\right)
>\mu,$ contradicting (\ref{cond}).

The same argument applies when $x_{k}<x_{l},$ completing the proof in this case.

Let now $q=2.$ If $t=4,$ then $\left\vert W\right\vert \geq e\left(
U,W\right)  \geq t=4,$ and so $W$ contains isolated edges, contradicting Claim
\ref{cl0}. Hence, $t=3,$ $\left\vert W\right\vert =3,$ and $G\left[  W\right]
$ is a path of order $3$. Now, the structure of $G$ is determined: $G$
consists of the graph $S_{m-4,3}$, a path $P$ of order $3,$ and a
$3$-matching, joining every vertex of $T$ to a separate triangle of
$S_{m-4,3}.$

At this point we apply again the above argument, completing the proof of
Theorem \ref{th1}.$\hfill\square$

\subsection*{\textbf{Concluding remarks}}

Theorem 3 in \cite{Nik07d} gives a result more general than just inequality
(\ref{in3}):

\begin{theorem}
Let $G$ be a graph of order $n$ with $\mu\left(  G\right)  =\mu$. If $G$ has
no $K_{2,k+1}$ for some $k\geq1,$ then
\[
\mu^{2}-\mu\leq t(n-1).
\]
Equality holds if and only if every two vertices of $G$ have exactly $k$
common neighbors.
\end{theorem}

This theorem is sharper than Theorem 3 in \cite{BaGu07}, and for some values
of $n$ and $k$ it is as good as one can get. However, in general, the maximal
$\mu\left(  G\right)  $ of $K_{2,k+1}$-free graphs $G$ of order $n$ is not
known at present.

Note that for $k>1,$ there may exist regular graphs with every two vertices
having exactly $k$ common neighbors: here is a small selection from
\cite{Roy94}:%

\[%
\begin{tabular}
[c]{||r|r|r||}\hline\hline
\emph{k} & \emph{n\ } & $\mu$\\\hline
$2$ & $16$ & $6$\\
$3$ & $45$ & $12$\\
$4$ & $96$ & $20$\\
$5$ & $175$ & $30$\\
$6$ & $36$ & $15$\\\hline\hline
\end{tabular}
\
\]
\bigskip

\textbf{Acknowledgement.} Thanks Laszlo Babai for the preprint \cite{BaGu07}%
.

\end{document}